\documentclass[a4paper,12pt,twoside]{article}
\usepackage{xcolor}
\usepackage[latin1]{inputenc}
\usepackage[T1]{fontenc}
\usepackage{amssymb}
\usepackage{dochier}
\usepackage[final]{lpretex}
\usepackage{title}
\usepackage{a4}
\usepackage{amscd}
\usepackage{caption}
\input xy
\xyoption{all}

\EndOfDump

\def\DHpartformat#1{(#1) }
\def\DHrefpart#1{(\DHRefpart{#1})}
\LBsetstyleof{partno}{arabic}

\let\define\def

\def\A {{\mathbb A}}  \def\C {{\mathbb C}}
  \def\F {{\mathbb F}}

  \def\P {{\mathbb P}} 
\def\Q {{\mathbb Q}}

\def\Z {{\mathbb Z}} 

\define \n {\mathbb N}
\define \z {\mathbb Z}
\define \q {\mathbb Q}
\define \PP {\mathbb P}

\def\sA {{\Cal A}}  
 \def\sE {{\Cal E}} \def\sF {{\Cal F}}

  \def\sO {{\Cal O}}
 
 \def\sT {{\Cal T}} 
  \def\sX {{\Cal X}}

\define \cN {\Cal N}
\define \cf {\Cal F}
\define \cg {\Cal G}
\define \cE {\Cal E}
\define \ce {\Cal E}
\define \cc {\Cal C}
\define \cV {\Cal V}
\define \cA {\Cal A}
\define \cK {\Cal K}
\define \cO {\Cal O}
\define \cF {\Cal F}
\define \cn {\Cal N}
\define \cI {\Cal I}
\define \sP {\Cal P}

\define \x {\xi}
\define \y {\eta}
\define \G {\Gamma}
\define \r {\rho}
\define \w {\omega}

\def \tZ {\widetilde Z}
\def\tX {\widetilde X}

\def \tC {\widetilde C}

\def \trho {\tilde {\rho}}

\def \tp {\widetilde{\mathbb P}}
\define \tH {\widetilde H}
\define \tG {\widetilde{\Gamma}}
\define \tW {\widetilde W}
\define \tF {\widetilde F}
\define \tm {\tilde m}
\define \St {\widetilde S}
\define \Xt {\widetilde X}
\define \tS {\widetilde S}
\define \tpsi {\tilde \psi}
\define \tL {\widetilde L}
\define \tE {\widetilde E}
\define \tl {\tilde l}
\define \tA {\widetilde A}
\define \tom {\tilde\omega}
\define \tT {\widetilde T}
\define \tB {\widetilde B}
\define \tf {\tilde f}
\define \tsA {\widetilde{\sA}}

\define \tM {\widetilde M}
\define \tphi {\widetilde{\phi}}
\define \trho {\widetilde{\rho}}
\define \tR {\widetilde R}
\define \tp {\tilde p}
\define \tq {\tilde q}
\define \tc {\tilde c}
\define \tsF {\widetilde {\sF}}
\define \tx {\tilde x}
\define \tg {\tilde g}
\define \tw {\tilde w}

\def\pd {\partial}

\def \Dx1 {\frac{\pd}{{\pd} x_1}}
\def \Dy1 {\frac{\pd}{{\pd} y_1}}
\def \Dz1 {\frac{\pd}{{\pd} z_1}}
\def \Dx2 {\frac{\pd}{{\pd} x_2}}
\def \Dy2 {\frac{\pd}{{\pd} y_2}}
\def \Dz2 {\frac{\pd}{{\pd} z_2}}

\def\q {\quad} 

\def\mapdiagr#1{\Big\searrow\rlap{$\raise 5pt\vbox{{\hbox{$\mkern -15mu\scriptstyle#1$}}}$}}   

\def\mapdiagl#1{\llap{$\raise 5pt\vbox{{\hbox{$\scriptstyle#1\mkern
-15mu$}}}$}\Big\swarrow}              

\def\Mapdiagr#1{\nearrow\rlap{$\lower 5pt\vbox{{\hbox{$\mkern
-15mu\scriptstyle#1$}}}$}} 

\def\Mapdiagl#1{\llap{$\lower 5pt\vbox{{\hbox{$\scriptstyle#1\mkern
-15mu$}}}$}\searrow} 

\def\Mapswr#1{\swarrow\rlap{$\lower 5pt\vbox{{\hbox{$\mkern
-15mu\scriptstyle#1$}}}$}}              

\def\Mapnwl#1{\nwarrow\rlap{$\lower 5pt\vbox{{\hbox{$\mkern
-15mu\scriptstyle#1$}}}$}}

\def \inj {\hookrightarrow}

\define \Rhook {\hookrightarrow}

\def \half {\raise1pt\hbox{$\scriptstyle
        \frac{1}{2}\displaystyle$}}

\def \x{{\sl X}\llap{$\mkern -2mu {\scriptstyle -}$}}
\def \Hom {\operatorname{Hom}}
\def \Spec {\operatorname{Spec}}

\def \Bl {\operatorname{Bl}}
\def \Pic {\operatorname{Pic}}
\def \Sing {\operatorname{Sing}}
\define \Kod {\operatorname{Kod}}
\define \dimension {\operatorname{dim}}
\define \codim {\operatorname{codim}}
\define \contr {\operatorname{contr}}
\define \rk {\operatorname{rank}}
\define \im {\operatorname{im}}
\define \Mor {\operatorname{Mor}}
\define \Cl {\operatorname{Cl}}
\define \Hilb {\operatorname{Hilb}}
\define \degree {\operatorname{deg}}
\define \mult {\operatorname{mult}}
\define \Aut {\operatorname{Aut}}
\define \NS {\operatorname{NS}}
\define \Gal {\operatorname{Gal}}
\define \ch {\operatorname{char}}
\define \Jac {\operatorname{Jac}}
\define \Km {\operatorname{Km}}
\define \Sec {\operatorname{Sec}}
\define \Stab {\operatorname{Stab}}
\define \Br {\operatorname{Br}}
\define \inv {\operatorname{inv}}
\define \tr {\operatorname{tr}}
\define \Frob {\operatorname{Frob}}
\define \Symn {\operatorname{Sym}^n}
\define \Ev {\sE^\vee}
\define \ordp {\operatorname{ord}_p}
\define \Supp {\operatorname{Supp}}
\define \Ann {\operatorname{Ann}}
\define \disc {\operatorname{disc}}
\define \Lie {\operatorname{Lie}}
\define \embdim {\operatorname{embdim}}

\def\Tr{\operatorname{Tr}}

\def\barK{\overline{K}}
\def\Tors{\operatorname{Tors}}

\def\hod#1#2#3#4{\ensuremath{\if#30 H^{#2}({#1},{\cal O}_{#1}) \else 
 H^{#2}(#1,\Omega^{#3}\if\relax{#4}\relax_{#1}\else _{#1/#4}\fi)\fi}}

\begin{document}
\title[9--nodal cubic threefolds]{The Hasse principle for 9--nodal cubic 3--folds}
\author{N. I. Shepherd-Barron}
\address{Dept. of Maths.,
  King's College, London WC2R 2LS,
 U.K.}
\email{Nicholas.Shepherd-Barron@kcl.ac.uk}
\maketitle
\begin{section}{Introduction}  
It has been known since the 19th century that if $X$ is
a cubic threefold over a field $K$ of characteristic zero
and $X\otimes\barK$ has just $d$
isolated singularities, then $d\le 10$.
Coray \cite{C} has shown that if $d$ is prime to 3,
then $X(K)$ is non--empty, while 
Colliot--Th{\'e}l{\`e}ne and Salberger \cite{CT--Sal} have shown that 
if $d=3$ and $K$ is a number field then 
$X$ satisfies the Hasse principle.
More recently Coray {\it{et al.}} \cite{C--L--SB--SD} have proved the
Hasse principle if $d=6$. This paper fills the gap, as far
as singular cubic threefolds are concerned, by showing that
if $d=9$, then (Theorem \ref{3.5} below)
the only obstruction to the Hasse principle
is the Brauer--Manin obstruction described in \cite{CT--San2}.
We do this by descent. More precisely, we show
first that, provided that $X$ has points everywhere locally,
universal torsors over the smooth locus $X^0$ of $X$
exist, that they are $K$--birational to cones over
certain singular cubic 7--folds and that all of them
satisfy the Hasse principle.

We also prove (Theorem \ref{7.5}) that on 10--nodal cubic threefolds
the only obstruction to weak approximation 
is the Brauer--Manin one, and prove a partial
such result (Proposition \ref{4.1}) in the
9--nodal case.


As mentioned, the proofs depend on the consideration
of various torsors, under tori.
What often makes various such torsors computable is
that the base variety is given by particularly simple 
(for example, linear) equations
inside some torus embedding, or rather an equivariant compactification
of some torsor under a torus. 9--nodal cubic threefolds follow
this pattern: they turn out to be hyperplane sections of some
Galois twist of the \emph{Perazzo cubic 4--fold} $P$, defined in
$\P^5$ by $x_1x_2x_3=y_1y_2y_3$. 
However, 10--nodal cubic 3--folds $X$ are different. They
are not, apparently, usefully embeddable in a toric variety,
and yet universal torsors over $X^0$ are simple; they are
birational to cones over the Grassmannian $G(2,6)$. 

\begin{acknowledgments}
I am grateful to Professors Browning, Colliot--Th{\'e}l{\`e}ne, Coray, Lewis, 
Skorobogatov and Swinnerton--Dyer for their valuable conversation, 
correspondence and encouragement, and to Professor Harpaz
for pointing out an error in an earlier version of this paper.

It is also a pleasure to
acknowledge the overwhelming
influence of Professor Manin on this area of mathematics.
\end{acknowledgments}
\end{section}
\begin{section}{Universal torsors over the smooth locus of
    a Perazzo cubic}\label{perazzo}
The object that renders 9--nodal cubic threefolds tractable
is the Perazzo cubic fourfold $P$ \cite{S--R}. This is given by the equation
$x_1x_2x_3=y_1y_2y_3.$
It contains nine 3--planes $L_{ij}$, given by
$x_i=y_j=0$. The complement $P^0=P-\cup L_{ij}$
is a torsor under a trivial $4$-dimensional torus $T$,
so $P$ is, by abuse of language,
a torus embedding. This makes it easy to compute universal torsors over it
and those of its twists that are also torus embeddings
(i.e., those that have $K$--points);
we shall do this explicitly in the next section.

Note that the nine 3--planes $L_{ij}$ are conjugate under the
wreath product
$\Gamma =S_3\wr S_2$; this is the subgroup of
the symmetric group $S_6$ generated by $S_3\times S_3$ and the involution 
$\iota$ of $S_3\times S_3$ that switches the two factors.
The first factor $S_3$ permutes the $x_i$, the second permutes
the $y_i$ and $\iota$ switches $x_i$ with $y_i$.
The singular locus $\Sing(P)$ consists of nine lines $l_{pq}$,
where $l_{pq}$ is given by $x_i=y_j=0$ for $i\ne p$,
$j\ne q$. The six points $w_1=(1,0,\ldots,0),\ldots,w_6=(0,\ldots,0,1)$
are $\Gamma$-conjugate; note that under the subgroup $S_3\times S_3$
of $\Gamma$ of index two they fall into two orbits.
Because
the complement $T$ of the 3--planes is a 4--dimensional torus and
the embedding $T\inj P$ is $T$--invariant, the class group 
$\Pic P^0 = \Cl(P\otimes\barK)$
is torsion--free of rank $9-\dim S =5$; it is generated by
the classes $L_{ij}$, subject to the relations of the form
$$0=(x_1/y_2)=L_{11}+L_{13}-L_{22}-L_{32}.$$

\begin{lemma}\label{1.1} The automorphism group scheme 
$G=\Aut (P,\sO(1))$ over $\Q$ is a split extension of $\Gamma$ by $T$.

\begin{proof}
The six points $w_i$ are distinguished as those points
where two or more of the $l_{pq}$ meet. So the connected
component $\Aut^0(P,\sO(1))$ preserves each of them,
and so preserves the 3--planes $L_{ij}$. Now the result is
obvious, with $T$ being the complement in $P$ of $\cup L_{ij}$.
\end{proof}
\end{lemma}

More generally, we define a \emph{Perazzo cubic fourfold}
to be a cubic $4$--fold $Y$ that is a $\Gal_K$--twist of $P$.
Its smooth locus will be denoted by $Y^0$.

\begin{proposition}\label{1.2} Every Perazzo cubic $Y$ satisfies the 
Hasse principle.

\begin{proof} Quadratic base extensions are harmless, 
so that we can assume that the 6 distinguished points fall
into two Galois orbits of three points each. Since it is known \cite{CT--Sal} 
that cubics with 3 conjugate nodes satisfy the Hasse principle, we are done.
\end{proof}
\end{proposition}

\begin{definition} A {\it{double--three}} is a configuration
of six 2--planes $L_i,M_j$ in $\P^8$, where $i,j=1,2,3$,
such that $\cup L_i$ and $\cup M_j$ each span $\P^8$ and
each intersection $L_i\cap(\cup M_j)$ and $M_j\cap(\cup L_i)$
consists of three non--collinear points.
\end{definition}

Note that any double--three has a unique decomposition into the union 
of two threes, where the 2--planes in each three are mutually disjoint.

\begin{proposition}\label{1.3} Suppose that $Y$ is a Perazzo cubic
with a $K$--point. 

\part[i] There exist universal torsors over $Y^0$.

\part[ii] Every such universal torsor
is $K$--birational to $\A^9_K$. 

\part[iii] The corresponding rational
map $\A^9_K-\to Y$ factors through the standard projection
$\A^9-\to\P^8_K$. 

\part[iv] The rational map $\P^8_K-\to Y\inj \P^5$ is given by a 
linear system of cubics passing doubly through a double--three.

\begin{proof} We know that $Y^0$ is a torus embedding $S_1\inj Y^0$, 
where $S_1$ is the complement of nine 3--planes in $Y$ and a torsor under the 
torus $S=\Aut^0(Y)$, and that $\barK[Y^0]^*=\barK^*$. Hence universal torsors 
over $Y^0$ exist and can be constructed according to the procedure
described in \cite{CT--San1}. That is, there is an exact sequence
$$1\to S_0\to M\to S\to 1$$
of tori, where $\hat S_0\cong \Pic(Y^0)$ as $\Gal_K$--modules and
$\hat M$ is the free module spanned by the classes of the 3--planes.
The coboundary map $S(K)\to H^1(K,S_0)$ is surjective
and the universal torsors over $Y^0$ form a torsor under
this $H^1$, so any universal torsor $\sT\to Y^0$ has the property
that $\sT|_S\to S$ is the pull--back of $\alpha :M\to S$
via the translation $\phi_x$ by $-x$, for some $x\in S(K)$. 

It therefore remains to describe the map $M\to S$.
We first do this in the untwisted case, where $Y=P$ and $S=T$.

The group $\Pic(P^0)$ is generated by the planes
$L_{ij}$ subject to the relations \hfill\break
$\sum_qL_{iq}-\sum_pL_{pj}
=(x_i/y_j)\sim 0$. Introduce nine new variables $z_{ij}$; then
there is a morphism $\pi:\A^9_K=\Spec K[\{z_{ij}\}]\to P$ defined  
by $x_i=\prod_qz_{iq}$ and $y_j=\prod_pz_{pj}$. Identify $M$
with the open subset of $\A^9$ given by $\prod z_{ij}\ne 0$;
then $\pi$ restricts to $\alpha$.

Consider the variables $z_{ij}$ as the entries of a $3\times 3$ matrix.
Then there are six triples of variables $z_{ab},z_{cd},z_{ef}$
such that no two of any triple lie in the same row or the same column.
For each such triple, let $L_{ab,cd,ef}$ be the 2--plane defined
by the vanishing of the other variables.
Then every cubic in the linear system defining
$\P^8-\to Y$ is double along each $L_{ab,cd,ef}$; it is clear 
that these 2--planes form a double--three, with 
$L_{11,22,33},L_{23,12,31}$ and $L_{21,32,13}$ forming one three.

For the general case, note that $S$ is the twist of $T$ by
some cocycle (in fact, homomorphism) $\Gal_K\to \Gamma$.
The action of $\Gamma$ on $T$ lifts to a $\pi$--equivariant linear
action on $\A^9$; taking the twists gives the results.
\end{proof}
\end{proposition}
\end{section}
\begin{section}{Geometry}
  In this section we investigate the basic geometry
  of 9--nodal cubic threefolds. Some of this material
  can also be found in \cite{C--T--Z}.

  We let $K$ denote a perfect field
of characteristic zero,
$\barK$ an algebraic closure of $K$ and $X$ a
cubic threefold over $K$ whose singular locus
$\Sing X$ consists
of exactly nine $\barK$-points.
We also assume that every Galois--conjugate subset
of $\Sing X$ has at least $2$ members, for else $X$ is $K$-rational
and there is nothing more to be said.
The smooth locus of $X$ will be denoted by $X^0$.

\begin{lemma}\label{2.1} 
\part[i] Every singularity of $X$ is an ordinary node. 

\part[ii] $X$ contains just nine 2--planes.

\part[iii] The 2--planes in $X$ and the points of $\Sing X$
form a $(9,4)$-configuration $G$ which is the 1--skeleton
of the cell decomposition of the 2--torus formed as the product of two
triangles.

\part[iv] The symmetry group of $G$ is $\Gamma$.

\begin{proof} For \DHrefpart{i}, we can assume that $K=\C$.

  Since cubic surfaces have at most $4$ isolated singularities,
  $X$ is not a cone and so,
  by I, Theorem 1.7 of \cite{Z}, the dual variety $X^\vee$
  of $X$ is a hypersurface. Moreover,
$$\deg X^\vee = 3.2^3-\sum_{v\in\Sing(X)}m(v),$$ 
where the class $m(v)$ of the singularity $(X,v)$ has the
property that $m(v)\ge 2$ and $m(v)=2$ if and only if $v$ is an 
ordinary node. Precisely, $m(v)=\mu(v)+\mu'(v)$, where $\mu$ is the 
Milnor number of an isolated hypersurface singularity $(X,v)$ and 
$\mu'$ is the Milnor number of a general hyperplane section \cite{T}. If 
$(X,v)$ is defined locally analytically (or formally)
in $\C^n$ by $f(x_1,\ldots,x_n)=0$, then
$$\mu(v)=\dim_\C\sO_{\C^n}/(\pd f/\pd x_1,\ldots,\pd f/\pd x_n).$$

If $\mu(v')\ge 2$ for some $v$ then $m(v)\ge 4$ for at least
two points $v$, so that
$$\deg X^\vee \le 24 -7\times 2 - 2\times 4=2.$$
But $\deg X^\vee\ge 3$, by the biduality theorem,
so that $\mu(v')=1$ for all $v$; this means that every
$v'$ is simple, of type $A_1$, and $v$ is
simple of type $A$.

Say that $X$ has $n_r$ points of type $A_r$,
so that $\sum n_r=9$.
Note that $m(A_r)=r+1$, so that
$$3\le 24 - \sum_{r\ge 1}n_r(r+1)=6-\sum_{r\ge 2} (r-1)n_r.$$
This gives $n_r=0$ for all $r\ge 5$ and $n_2+2n_3+3n_4\le 3$,
so that $n_1\ge 6$. 

Choose an $A_1$-point $P$. Projecting from $P$
identifies $\Bl_PX$ with $\Bl_C\P^3$, where
$C$ is a reduced curve of bidegree $(3,3)$
on a smooth quadric with $n_1-1$ points of type
$A_1$ and $n_r$ of type $n_r$ for all $r\ge 2$.
Suppose that $\tC=\sqcup \tC_j\to C$ is the normalization;
then
$$\sum\chi(\sO_{\tC_j})=\chi(\sO_C)+(n_1-1).1+n_2.1+n_3.2+n_4.2
=5+n_3+n_4.$$
So, if $n_3+n_4>0$, then $C$ is the transverse union of six lines,
and then $n_3=n_4=0$, contradiction. Therefore $n_3=n_4=0$,
so that $C=\sum_1^5 C_j$ where $C_5$ is a conic and the other
components are lines.
Then also $n_2=0$ and \DHrefpart{i} is proved.
\begin{remark}
If there is a tenth singular point, then 
it is known classically that $X$ is a form of the Segre cubic.
\end{remark}
For \DHrefpart{ii}, note that the description
of the curve $C$ in the proof of \DHrefpart{i}
shows that there are exactly four 2--planes on $X$ through
each node.
Each plane in $X$ contains 4 nodes of $X$, so that
there are just 9 planes in $X$.

\DHrefpart{iii}: From the description of $C$ it also follows that
through each node of $X$ there are 4 lines in $X$ that pass through 
a further node and every such line is the intersection of a unique pair 
of 2--planes in $X$. Each of the 4 nodes in a 2--plane $L$ lies on two
lines of the form $L\cap M$, since each line in $C$ meets
two others. This proves \DHrefpart{iii}.

\DHrefpart{iv}: it is clear that the symmetry group of $G$ is
isomorphic to $S_3\wr S_2$.
\end{proof}
\end{lemma}

\begin{remark} We can also view $G$ as a graph whose vertices
  are the planes in $X$ and whose edges are the pairs of planes that
  meet in a point.
\end{remark}

\begin{proposition}\label{2.2} Any $9$--nodal cubic $3$--fold $X$
over $K$ is $K$--isomorphic to a hyperplane section of some Perazzo 4--fold.

\begin{proof} Assume first $K=\barK$. By \ref{2.1} there are nine
2--planes in $X$; denote them by $L_{ij}$, where
$i,j\in \Z/3$ and $L_{ij}$ meets $L_{i\pm 1,j}$ and
$L_{i,j\pm 1}$, each in a line. So for all $i$ and $j$,
$\sum_iL_{ij}$ and $\sum_jL_{ij}$
are hyperplane sections. Say $\sum_jL_{ij}=(x_i=0)$
and $\sum_iL_{ij}=(y_j=0)$. Then $\sum_{ij}L_{ij}$ is the
complete intersection $\prod x_i=\prod y_j=0$, so that
$X$ is a member of the pencil $|3H-\sum_{ij}L_{ij}|$.
Hence the equation of $X$ is $\alpha \prod x_i=\beta\prod y_j$,
as required.

Next, we check (still assuming $K=\barK$) that if $X,X'$
are isomorphic sections of $P$, then there is an automorphism
of $P$ taking $X$ to $X'$.

There exists $\sigma\in PGL_6$ such that $\sigma(X)=X'$.
Put $P'=\sigma(P)$; then there is a hyperplane $H$ such that
$P.H=P'.H$, and it is enough to find $\tau\in PGL_6$
such that $\tau(P)=P'$ and $\tau|_H=1$. Put 
$G=\{\tau\in PGL_6|\tau|_H=1\}$; then $G$ acts transitively on
$\P^5-H$. 

Both $P$ and $P'$ contain nine 3--planes $M_{ij},M'_{ij}$,
respectively; we can take them to be ordered so that
$M_{ij}.H=M'_{ij}.H$ for all $i,j$. Suppose that $P$, resp. $P'$,
is given by $F=0$, resp. $F'=0$, where $F=x_1x_2x_3-y_1y_2y_3$; then 
we can take $M_{ij}$ to be given by $x_i=y_j=0$.
Put $v=(1,0,0,0,0,0)$, so that $v=\cap_{i\ne 1}L_{ij}$.
So, after applying a suitable $\tau\in G$, we have
$\cap_{i\ne 1}L_{ij}=\cap_{i\ne 1}L'_{ij}$. Since
$L_{ij}.H=L'_{ij}.H$, we now have
$L_{ij}=L'_{ij}$ for all $i\ne 1$.
So $F'\in \cap_{i\ne 1}(x_i,y_j)=(x_2x_3,y_1y_2y_3)$; the equality
of these two ideals is verified by noting that both define
subschemes of codimension 2 and degree 6
and that obviously
$(x_2x_3,y_1y_2y_3)\subseteq\cap_{i\ne 1}(x_i,y_j)$.
Hence $F'=lx_2x_3+\alpha y_1y_2y_3$, where $\alpha\in K$
and $l$ is linear. Say that $H$ is defined by $m=0$; it is
then easy to see, since $X$ is irreducible, that, after
rescaling, $m=l-x_1$ and $F'=(x_1+m)x_2x_3-y_1y_2y_3$.
Then $\tau\in GL_6$ given by
$$x_1\mapsto x_1+m, x_i\mapsto x_i{\quad{\textrm{for}}}\quad i\ne 1,
y_j\mapsto y_j$$
has the required effect.

Now drop the assumption that $K=\barK$. There is an embedding 
$\phi:X_{\barK}\inj P_{\barK}$ with the property that for all
$\sigma\in\Gal_K$, there exists 
$\psi_\sigma\in \Aut P(\barK)$ such that
$\phi^\sigma=\psi_\sigma\circ\phi$. Then
$\psi_{\sigma\tau}\phi =(\psi_\sigma)^\tau\psi_\tau\phi$
for $\sigma,\tau\in\Gal_K$.

Write $\omega=\psi_{\sigma\tau}^{-1}(\psi_\sigma)^\tau\psi_\tau$.
So $\omega\phi=\phi$. That is, $\omega$ is an automorphism
of $P_{\barK}$ that acts trivially on the hyperplane $H$
cutting out $X_{\barK}$. 

We want to show that any such $\omega$ is the identity. For this,
we can assume once more that $K=\barK$ and that $\omega\ne 1$. 
After moving $H$ by an element of $S$
we can suppose also that $(1,\ldots,1)\in H$. Then $\omega\in\Gamma_0$,
where $\Gamma_0\subseteq \Aut P$ is a copy of $\Gamma$ that splits the
surjection $\Aut P\to \Gamma$.
Since $\omega$ is trivial on $H$, it is conjugate to
$(s,1)\in S_3\times S_3\subseteq\Gamma_0$, where $s$ is a transposition;
we can take $s=(12)$. Then $H$ is given by $x_1=x_2$,
which contradicts the fact that $P.H$ has isolated singularities.
So $\omega =1$ and
$(\psi_\sigma)\in Z^1(\Gal_K,\Aut P(\barK))$.
It is now easy to see that 
$X$ embeds into the twist of $P$ by $(\psi_\sigma)$.
\end{proof}
\end{proposition}

Now suppose that $X$ is a hyperplane section of $Y$ and that
$\Sing X=X\cap\Sing Y$. So $X$ is a 9--nodal cubic threefold.
Let $\tX\to X$ the blow--up of the nodes.
It is then easy to see that the Betti numbers of $\tX$
are determined by $e(\tX)=-6+9.4=30$ and $b_3=0$,
so that $\rk \Cl(X\otimes\barK)=5$. 

\begin{lemma}\label{2.3}
  \part[i] The natural map
  $\Cl(Y_{\barK})\to\Cl(X_{\barK})$ is a
  $\Gal_K$--isomorphism.
  \part[ii] $\Cl(X_{\barK})$ is generated by the classes $L_{ij}$
  subject to the relations $R_i-C_j=0$, where $R_i=\sum_i L_{ij}$
  and $C_j=\sum_j L_{ij}$.

\begin{proof} We check first that the 2--planes on $X$ generate
$\Cl(X)$.

Fix a node $P$ on $X$; then, as before, $\Bl_PX$
is identified with $\Bl_C\P^3$ 
and we see that
$\Cl(\Bl_PX)$ is generated
by $L_1,\ldots,L_4,Q,H$, where $L_i$ is the strict transform
of a plane that projects to a line $C_i$ in $C$, $Q$ the strict transform
of a quadric cone projecting to the conic $C_5$ in $C$
and $H$ is the pull--back
of the hyperplane class on $\P^3$. 

Let $H_1$ be the hyperplane class on $X$. Then $H_1-E\sim H$
in $\Cl(\Bl_PX)$
and there are 2--planes $L',L''$ on $X$ such that
$H_1\sim L_1+L_3+L'$ and $H_1\sim Q+L''$ in $\Cl(X)$.
Since $\Cl(X)\cong \Cl(\Bl_PX)/\Z.[E]$, we get relations
$H\sim H_1\sim L_1+L_3+L'$ and
$Q\sim H_1-L''$ in $\Cl(X)$, and the planes on $X$ do
generate $\Cl(X)$. The relations $R_i=C_j$ follow from the
observation that $R_i\sim H_1\sim C_j$.

In $Y$, the 3--planes form the boundary $Y-U$, as described
above, and so generate $\Cl(Y)$. Hence $\Cl(Y)\to \Cl(X)$
is surjective.
Since both groups have rank 5, it is enough to prove that $\Cl(X)$
is torsion--free. Since $X$ has isolated hypersurface singularities
and is 3--dimensional, $\pi_1(X^0)\to\pi_1(X)$ is an isomorphism.
Since $\pi_1(X)=1$, we are done.
\end{proof}
\end{lemma}
\begin{lemma} The image $W$ of $\Gal_K$ on $\Cl(X_{\barK})$
  equals its image in $\G$.
  \begin{proof} This is a consequence of the fact, which has
    already been remarked, that the nine 2--planes in $X$
    generate $\Cl(X_{\barK})$.
  \end{proof}
\end{lemma}
\end{section}
\begin{section}{The Hasse principle}
  \label{Hasse}
Now suppose that $K$ is a number field and that $X$ is a $9$--nodal 
cubic threefold over $K$. We know that $X$ is a hyperplane 
section of some Perazzo cubic $4$-fold $Y$ over $K$.

\begin{lemma}\label{3.1} Suppose that $\sT\to Y^0$ is a universal torsor.
Then so is $\sT\times_{Y^0}X^0\to X^0$. Moreover, every universal
torsor on $X^0$ arises in this way.

\begin{proof} Let $S_0$ denote the torus whose character group 
$\hat S_0$ is $\Pic(Y^0)$. There is a commutative diagram (cf. \cite{CT--San2})
$$\xymatrix{
{0}\ar[r] & {H^1(K,S_0)} \ar[r]\ar[d]_{=} & {H^1(Y^0,S_0)}\ar[r]^-{\chi}\ar[d] &
{\Hom_{\Gal_K}(\hat S_0,\Pic Y^0)}\ar[d]^{\cong}\\
{0}\ar[r] & {H^1(K,S_0)}\ar[r] & {H^1(X^0,S_0)}\ar[r]^-{\chi} & 
{\Hom_{\Gal_K}(\hat S_0,\Pic X^0).}
}$$ 
By definition, a torsor under $S_0$ is universal if the image of
its class under $\chi$ is the identity, and now the result is immediate.
\end{proof}
\end{lemma}

{\it{Now assume that $X$ has $K_v$--points
for all places $v$ of $K$.}}
That is, $X(\A_K)$ is not empty, where $\A_K$
is the ring of ad\`eles of $K$.

\begin{proposition}\label{3.2} There is a universal torsor over $X^0$.

\begin{proof} We know that $X$ is a section of $Y$; since $Y$ satisfies
the Hasse principle it has a $K$--point, and so there is 
a universal torsor over $Y^0$. Now use Lemma \ref{3.1}.
\end{proof}
\end{proposition}

\begin{proposition}\label{3.3} Every universal torsor over $X^0$ is 
$K$--birational to the cone over a cubic 7--fold that is singular along 
a double--three.

\begin{proof} Immediate from Proposition \ref{1.3}.
\end{proof}
\end{proposition}

\begin{proposition}\label{3.4} The universal torsors $\sT$ over $X^0$ 
satisfy the Hasse principle.

\begin{proof} Say that $\sT$ is $K$--birational to the cone
over the cubic 7--fold $Z$. Since quadratic extensions are
harmless, we can suppose that the given double--three along which
$Z$ is singular splits into two threes. Then $Z$ has three conjugate 
singular points, and so \cite{CT--Sal} satisfies the Hasse principle.
\end{proof}
\end{proposition}

For any $K$-variety $V$, there is a pairing
$$\Br(V)\times \prod V(\A_K)\to\Q/\Z$$
given by
$$(\sA,(P_v))\mapsto \sum_v\inv_v(\sA(P_v)).$$
We denote by $V(\A_K)^{\Br}$ the subset of $V(\A_K)$ that is the kernel
of this pairing. The set $V(K)$ lies naturally in $V(\A_K)^{\Br}$.
If the non-emptiness of $V(\A_K)^{\Br}$ implies
that of $V(K)$ then we say that
``the only obstruction to the Hasse principle on $V$
is the Brauer--Manin obstruction''.

\begin{theorem}\label{3.5} 
The Brauer--Manin obstruction to the Hasse principle on $X$ is the only one.

\begin{proof} By definition \cite{CT--San2}, the Brauer--Manin
obstruction to the existence of a $K$--point on $X$ (equivalently,
since $X$ is a cubic, on $X^0$) is the existence, for all
$(x_v)\in X^0(\A_K)$, of an element $A$ of $\Br(X^0)$
such that $\sum_v\inv_v(A)\ne 0$. Assume that this obstruction
vanishes; then the proof of Th{\'e}or{\`e}me 3.8.1
of \emph{loc. cit.} shows that there is a universal torsor over $X^0$ 
with a point everywhere locally. By Proposition \ref{3.4}, we are done.
\end{proof}
\end{theorem}

Colliot--Th{\'e}l{\`e}ne points out that the following variant of Th. 3.8.1 
of \emph{loc. cit.} is valid, where $\tX$ denotes a smooth compactification 
of $X^0$ (for example, the blow--up of the nodes of $X$). 

\begin{lemma}\label{3.7} Suppose that $Z$ is a projective variety 
over $K$ with only nodes and that $\dim Z\ge 2$. Denote by
$Z^0$ its smooth locus and $\tZ\to Z$ the blow--up of the nodes. 
Then the natural map $\Br(\tZ)\to\Br(Z^0)$ is an isomorphism.

\begin{proof} We first prove this when $K=\barK$.

  Brauer groups are torsion, so it is enough to
prove this for the $n$--torsion subgroups.
The Kummer sequence shows that then it is enough
to prove the surjectivity of $H^2(\tZ,\mu_n)\to H^2(Z^0,\mu_n)$.
Let $E=\sum E_i$ be the exceptional locus in $\tZ$ and
$j:Z^0\to\tZ$ the inclusion. Then the
purity theorem shows that $j_*\mu_n=(\mu_n)_{\tZ}$,
$R^1j_*\mu_n$ is locally isomorphic to $(\mu_n)_E$ and
that $R^qj_*\mu_n=0$ for $q\ge 2$. Since each $E_i$ is simply
connected, it follows that $H^1(\tZ,R^1j_*\mu_n)=0$. Now 
the Leray spectral sequence 
$E_2^{pq}=H^p(\tZ,R^qj_*\mu_n)\Rightarrow H^{p+q}(Z^0,\mu_n)$ 
gives the result.

The general case then follows from the facts that
$\Br(V)/\Br(K)$ is naturally isomorphic to
$H^1(\Gal_K,\Pic(V_{\barK}))$ and that 
the homomorphism $\Br(\tZ_{\barK})\to\Br(Z^0_{\barK})$
is $\Gal_K$--equivariant.
\end{proof}
\end{lemma}

\begin{proposition}\label{3.6} (Colliot-Th{\'e}l{\`e}ne.)
If there is no Brauer--Manin obstruction
using $\Br(\tX)$ to the existence of a $K$--point on $X$,
then there is none using $\Br(X^0)$.

\begin{proof} Choose a finite set $\sA_1,\ldots,\sA_n$ of 
representatives of the elements of
the finite group $\Br(X^0)/\Br(K)$.  Over an open subscheme
$\Spec(\sO)$ of the spectrum of the ring of integers of $K$, the variety $X^0/K$ has
a model ${\sX^0}/\sO$ such that  ${\sX^0}(\sO_v) \neq \emptyset$
for each $v \in \Spec(\sO)$, and such that each $\sA_i$ extends
to an element of ${\Br}({\sX^0})$.
Using Harari's ``formal lemma'' (Lemme 2.6.1 of \cite{H}, but 
see also p. 225 of \cite{CT--San2}) and the hypothesis that there is no
Brauer-Manin obstruction on $\tX$, we find a finite set $S$ of places, which we
may assume contains all places not in $\Spec(\sO)$,
and local points $M_v \in X^0(K_v)$ for $v \in S$,
such that $$ \sum_{v \in S} {\inv_v}(\sA_i(M_v))=0$$
for each $i =1, \ldots, n$.

Now pick any set of integral points $M_v \in {\sX^0}(\sO_v)$
for $v \notin S$.  Then for each $i \in \{1, \cdots, n \}$, the sum
$\sum_{v} {\inv_v}(\sA_i(M_v))$ vanishes (the sum is over all places of $K$),
thus completing the proof.
\end{proof}
\end{proposition}
\end{section}
\begin{section}{Weak approximation}
  Recall that a 
variety $V$ over $K$ satisfies \emph{weak approximation} (WA) if
for every finite set $S$ of places of $K$,
$V(K)$ is dense in $\prod_{v\in S}V(K_v)$.
If $V$ is complete, this is equivalent to $V(K)$
being dense in $V(\A_K)$. Moreover,
we say that the Brauer--Manin obstruction
to WA on $V$ is the only one if $V(K)$ is dense in $V(\A_K)^{\Br}$.
This is equivalent to the density of $V(K)$ in the image of
$V(\A_K)^{\Br}$ under every projection 
$V(\A_K)\to \prod_{v\in S}V(K_v)$, for every finite set $S$.

Now assume that $X$ and $K$ are as in Section \ref{Hasse}.

\begin{lemma} $X^0(\A_K)^{\Br}$ is dense in $X(\A_K)^{\Br}$.
\begin{proof} This is Corollary 1.2 of \cite{CT--Sk}.
\end{proof}
\end{lemma}

\begin{proposition}\label{4.1} Assume that
$W$ lies in the index 2 subgroup $S_3\times S_3$
of $\Gamma$. Then the only
obstruction to weak approximation on $X$ is the Brauer--Manin one.
\begin{proof} 
Take a point $(P_v)\in X(\A_K)^{\Br}$. By the lemma,
we can approximate $(P_v)$ by $(M_v)\in X^0(\A_K)^{\Br}$.
By the version of descent given in Proposition 1.3 of
\cite{CT--Sk}, there is a universal torsor $\pi:T\to X^0$
such that $(M_v)$ lifts to a point $(\tM_v)\in T(\A_K)$. We know that $T$ is
$K$--birational to the cone over a cubic 7--fold $Z$, and
that $Z$ is singular along a double--three. The hypothesis on
the Galois group means that the double--three splits into 
two threes, so that $Z$ has three conjugate singularities.
Since this gives WA for $Z$, by \cite{CT--Sal}, $T$ then satisfies WA,
so that there exists $t\in T(K)$ close to $(\tM_v)$.
Then $x:=\pi(t)\in X^0(K)$ is close to $(M_v)$,
and then $x$ is close to $(P_v)$.
\end{proof}
\end{proposition}
\end{section}
\begin{section}{Computing the Brauer group and obstructions}
\medskip
Suppose that $X$ is a $9$--nodal cubic $3$--fold over $K$
and $\tX$ its blow-up at the nodes. Put 
$\Pic(X^0\otimes\barK)= P$.  Then, by Lemma \ref{3.7} and the 
Hochschild--Serre spectral sequence, $\Br(X^0)/\Br(K)$ is isomorphic 
to $H^1(\Gal_K,P)=H^1(W,P)$. 

\begin{proposition}\label{5.1} 
$\Br(X^0)/\Br(K)$ is either trivial or of order $3$.

\begin{proof} 
Pick a general projection $\pi:X\to\P^1_K$, so that the generic
fibre $X_\eta$ is a smooth cubic surface
and all geometric fibres are irreducible. 
Put $K(\P^1)=L$. The main result of \cite{SD1} is that $\Br(X_\eta)/\Br(L)$ 
is a subgroup of either $(\Z/2\Z)^2$ or $(\Z/3\Z)^2$.

We show first that the natural homomorphism
$$\phi:\Br(X^0)/\Br(K)\to\Br(X_\eta)/\Br(L)$$ is injective.
For this, suppose that $\alpha\in\Br(X^0)$ restricts
to $\pi^*\beta$ in $\Br(X_\eta)$. For any closed point $M$
of $\P^1_K$, the image of $\alpha$ under the residue
map $\pd_{\pi^{-1}(M)}:\Br(K(X))\to H^1(K(\pi^{-1}(M)),\Q/\Z)$ is zero.
Since $K(M)$ is algebraically closed in $K(\pi^{-1}(M))$,
it follows that $\pd_M(\beta)=0$, so that 
$\beta\in\Br(\P^1_K)=\Br(K)$. So $\phi$ is injective.

Suppose next that there is 2--torsion in
$\Br(X^0)/\Br(K)$. Then the same is true after any
base extension of odd degree, so that we can assume that 
$X$ has a $K$--rational node. Then $X$ is $K$-rational,
so that, by Lemma \ref{3.7}, $\Br(X^0)/\Br(K)=0$,
contradiction.
So $\Br(X^0)/\Br(K)$ has odd order.

Since $\Gamma\cong S_3\wr S_2$,
we can assume that $W$ is non-trivial and
lies in the unique Sylow $3$-subgroup of $\Gamma$, which is 
$A_3\times A_3$. So, up to $\Gamma$--conjugacy, there are three
possibilities for $W$: $A_3\times 1$; the diagonal copy 
$\Delta$ of $A_3$ in $A_3\times A_3$; and $A_3\times A_3$.

Recall that $P$ is generated by the $L_{ij}$ subject
to the relations $R_i-C_j=0$. We consider the three
possibilities separately.

\begin{enumerate}
  \item $W=A_3\times 1$. 
Put $P_1=\oplus_{j\ne 3}\Z.L_{ij}$.
This is a permutation $W$--module and
there is a short exact sequence
$$0\to\Z.(C_1-C_2)\to P_1\to P\to 0.$$
Since $H^2(W,\Z)\cong\Z/3\Z$ and $H^i(W,P_1)=0$ for $i=1,2$, 
it follows that $H^1(W,P)$ is isomorphic to $\Z/3\Z$.

\item $W=\Delta$. Put $T=s-1$, where $s=((123),(123))$ is a 
generator of $\Delta$, and $N=\sum s^i$. So, for any $W$--module $B$,
$H^1(W,B)\cong \ker N_B/\im T_B\cong\Tors(\coker T_B)$.
Let $F$ denote the free $\Z$-module on the $L_{ij}$,
so that there is a short exact sequence
$$0\to Q\to F \to P\to 0$$
of $W$-modules, which defines the submodule $Q$ of $F$.
Inspection shows that for $i\ne j$
$L_{ii}-L_{jj}=\pm T(L_{ii})$ or $\pm T(L_{jj})$ and
$L_{ik}-L_{kj}=\pm T(L_{ik})$ or $\pm T(L_{kj})$.
Hence $Q\subseteq \im(T_F)$, so that
$\coker T_P\cong F/(Q+\im T_F)\cong\coker T_F$.
Since $H^1(W,F)$ vanishes, so does $H^1(W,P)$.

\item $W=A_3\times A_3$. By Lemma 5 of \cite{SD1}
$\Br(X_\eta)/\Br(L)$ is then of order at most 3, so the same is
true of $\Br(X^0)/\Br(K)$.
\end{enumerate}
\end{proof}
\end{proposition}

\begin{proposition}\label{5.2} Suppose that $W\subseteq A_3\times
  A_3$. Then $\Br(X^0)/\Br(K)\cong \Z/3\Z$ if $W$ is $\Gamma$--conjugate
to either $A_3\times A_3$ or $A_3\times 1$, and is trivial
otherwise. Moreover, when $\Br(X^0)/\Br(K)$ is non--trivial
and $K$ contains a cube root of unity, there is an explicit
description (given in the course of the proof)
of a non--trivial element.

\begin{proof} We use the results and notation of \cite{SD2},
especially Lemma 2 of \emph{loc. cit.}

Assume that $\Br(X^0)/\Br(K)$
is non--trivial. There is a cubic extension $K_1/K$
over which the divisors $L_{11},L_{22},L_{33}$ are defined.
The image $W'$ of $\Gal_{K_1}$ in $\G$
is of order $3$ and is generated by $\sigma=((123),1)$.
Choose a hyperplane section $H$ defined over $K$ and
put $D=\sum L_{ii} -H$. Then $\sum_j\sigma^j(D)$
is principal; say $\sum\sigma^j(D)=(f)$, with $f\in K(X)$.
Then, according to \emph{loc. cit.}, there is a non--trivial
element $\sA$ of $\Br(X^0)$ such that, for every adelic point
$(P_v)$ on $X^0$, with no $P_v$ in the support of $(f)$,
$\sum_v\inv_v(\sA(P_v))=\sum_v(f(P_v),K_1/K)_v$, where the
summands on the right are the norm residue symbols.
\end{proof}
\end{proposition}
%
%
%
\end{section}
\begin{section}{An aside: weak approximation for 10 nodes}
  Coray has shown \cite{C} that a 10-nodal cubic threefold $X$ has
  $K$-points. We show here that also weak approximation holds.
  
We start by recalling a construction from \cite{D--O}. There the ground field is $\C$, 
but this part of \emph{loc. cit.} is valid over any field $K$, or even over $\Z$.

Regard $T=\mathbb G_{m,K}^6$ as the group of $6\times 6$ diagonal
matrices acting in the obvious way on the 6--dimensional $K$--vector space $V$.
Let $U$ be the standard 2--dimensional representation of $SL_2$ over $K$.
Then $G=(SL_2\times T)/\mu_2$ acts on the space
$M=U\otimes V\cong\mathbb A^{12}$ of $2\times 6$ matrices, where
$\mu_2$ is embedded diagonally. Regard $SL_2$ as acting on the left and $T$ on the right.
Let $M^0$ be the locus of matrices of rank 2; then
$SL_2$ acts freely on $M^0$ and there is a geometric quotient
$SL_2\backslash M^0$ isomorphic to the punctured cone 
$\widetilde{Gr}$ over the Grassmannian $Gr(2,V)$ in its Pl\"ucker embedding.
If $M^{00}$ is the locus of matrices where no column
is zero and no 3 columns are proportional, then 
$G$ acts freely on $M^{00}$. There is a geometric
quotient $M^{00}/G$, which is isomorphic to the smooth locus
$\Sigma^0$ of the Segre cubic 3-fold $\Sigma$ given by $\sum_1^6x_i^3=\sum_1^6x_i=0$.
We can identify $\Sigma^0=M^{00}/G$ with $(SL_2\backslash M^{00})/S$, 
where $S=T/\mu_2$ with $\mu_2$ embedded in the diagonal copy of
$\mathbb G_m$ in $T$. It follows that the $S$-torsor $\widetilde{Gr}\to\Sigma^0$
is a universal torsor over $\Sigma^0$.
\def\barQ{\overline{\mathbb Q}}
The basic geometry of $\Sigma$ is described in \cite{S--R}, p.169.
It is the image of $\P^3$ under the rational map defined by
the linear system of quadrics that pass through five given
points in general position. It has 10 nodes and 15 planes, all defined over $\Q$,
and over $\barQ$ it is the unique 10-nodal cubic threefold. That is, 
every 10-nodal cubic threefold $X$ over $K$ is a twist of $\Sigma$.

\begin{lemma}\label{7.1}

\part[i] $\Aut (\Sigma)= \Aut (\Sigma\otimes \bar\Q)\cong S_6$.

\part[ii] Every $1$-cocycle $\Gal_Q\to\Aut (\Sigma)$ is a homomorphism.

\part[iii] Every $10$-nodal cubic threefold $X$
over $K$ is the twist of $\Sigma$ by a homomorphism $\Gal_K\to S_6$.

\part[iv] Given such an $X$ over $K$, there is a separable
sextic $K$-algebra $L$ such that $X$ is defined in $\P(L)$
by the equation $\Tr(z)=\Tr(z^3)=0$.

\begin{proof} \DHrefpart{i} 
The projective dual of $\Sigma\otimes\barQ$ is a quartic threefold $T$
that is the Satake compactification of the moduli space of principally
polarized Abelian surfaces with level 2 structure. The Satake boundary
is $\Sing T$, which is a $(15_3,15_3)$ configuration of points and lines.
The automorphism group of this is clearly $Sp_4(\F_2)\cong S_6$,
so that there is a homomorphism $f: \Aut (\Sigma\otimes\barQ)\to S_6$.
Restrict $f$ to the copy of $S_6$ in
$\Aut(\Sigma)\subseteq\Aut(\Sigma\otimes\barQ)$; this is clearly an
isomorphism, so that $f$ is split. Since the lines in $\Sing T$ are the
images of the planes in $\Sigma$, any $s\in \ker f$ preserves every
plane in $\Sigma$. Since each node of $\Sigma$ is (in many ways) the
intersection of two of these planes, $s$ fixes every node. It is then
clear that $s=1$.

\DHrefpart{ii} This follows at once from Lemma 2.1 and the definition
of a 1-cocycle.

\DHrefpart{iii} and \DHrefpart{iv} are now immediate. 
\end{proof}
\end{lemma}

Now suppose that $X/K$ is a 10--nodal cubic 3--fold with smooth
locus $X^0$. We have just seen that $X$ is $K$--isomorphic to a 
twist $\Sigma_\rho$ of $\Sigma$ by a homomorphism $\rho:\Gal_K\to S_6$.
Then the construction above can be twisted to show that
$X^0$ is the geometric quotient $M_\rho^{00}/G_\rho$.

\begin{proposition} There is a universal torsor $T_0$ over $X^0$ whose
total space is isomorphic to $\widetilde{Gr}$.
\begin{proof} The morphism $\widetilde{Gr}\to \Sigma^0$ can be twisted 
by $\rho$. The twist of $Gr$ by $\rho$ is $Gr(2,V_\rho)$, which is
isomorphic to $Gr$, by Hilbert's Theorem $90$. It is also possible 
to twist the action of the torus $S$ on $\widetilde{Gr}$ by $\rho$,
and we are done, since $\widetilde{Gr}$ has trivial Picard group.
\end{proof}
\end{proposition}

\begin{corollary} Every universal torsor $T$ over $X^0$ is $K$-birational
to a line bundle over a twist of $Gr$.
\begin{proof} Given $T_0$, the other universal torsors over $X^0$
are classified by $H^1(K,S)$. Since $S$ acts by right multiplication
on $Gr$ and the Pl\"ucker line bundle $\sO(1)$ over
$Gr$ is $S$-linearized, the result follows.
\end{proof}
\end{corollary}

\begin{corollary} Both the Hasse principle and weak approximation
hold for every universal torsor $T$ over $X^0$.
\begin{proof} Every twist $V$ of $Gr$ is homogeneous under its
automorphism group, and the point stabilizers are connected.
Hence the Hasse principle and weak approximation
hold for $V$, and then for line bundles over it.
\end{proof}
\end{corollary}

\begin{theorem}\label{7.5} The only obstruction to weak approximation
on $X$ is the Brauer--Manin one.
\begin{proof} Exactly as for the class of 9--nodal cubics
that we handled earlier, in \ref{4.1}.
\end{proof}
\end{theorem} 

\end{section}

\end{document}